\documentclass{article}

\usepackage{amsmath, amsthm, amssymb, tikz, hyperref, enumitem,array,mathtools}
\usepackage{amsfonts}
\usepackage{fullpage}
\usepackage{shuffle}

\newtheorem{theorem}{Theorem}[section]
\newtheorem{proposition}[theorem]{Proposition}
\newtheorem{lemma}[theorem]{Lemma}
\newtheorem{corollary}[theorem]{Corollary}

\newtheorem{remark}[theorem]{Remark}
\newtheorem{observation}[theorem]{Observation}

\theoremstyle{definition}
\newtheorem{example}[theorem]{Example}  
\newtheorem{definition}[theorem]{Definition}   

\DeclareMathOperator{\Des}{Des}
\DeclareMathOperator{\des}{des}
\DeclareMathOperator{\ev}{ev}
\DeclareMathOperator{\st}{st}
\DeclareMathOperator{\maj}{maj}
\DeclareMathOperator{\comaj}{comaj}

\newcommand{\perm}{p}
\renewcommand{\S}{\mathcal{S}}
\newcommand{\Co}{\mathcal{C}}
\newcommand{\bij}{\phi}
\newcommand{\bijS}{\phi_S}

\newcommand{\id}{\iota}
\newcommand{\card}[1]{\left|#1\right|}

\newcommand{\co}{\Des}

\def\bl{\textcolor{blue}}
\def\rd{\textcolor{red}}

\DeclareMathOperator{\Qsym}{Qsym}

\renewcommand{\a}{\mathtt{a}}
\renewcommand{\b}{\mathtt{b}}
\renewcommand{\c}{\mathtt{c}}
\renewcommand{\d}{\mathtt{d}}
\newcommand{\e}{\mathtt{e}}
\newcommand{\f}{\mathtt{f}}

\def\NN{{\mathbb N}}
\def\ZZ{{\mathbb Z}}
\def\QQ{{\mathbb Q}}

\newcommand{\shc}[2]{#1 \shuffle_c #2}

\newcommand{\delete}[1]{}

\title{Shuffles of permutations in cycle form}
\author{
Ron M. Adin~\thanks{Department of Mathematics, Bar-Ilan University, Ramat-Gan 52900, Israel. 
E-mail: {\tt radin@math.biu.ac.il}.
Partially supported by Dartmouth-Kalaniyot in the Sciences (DKS) Grant \#1013.}
\and
Sergi Elizalde~\thanks{Department of Mathematics, Dartmouth College, Hanover, NH 03755, USA. 
E-mail: {\tt sergi.elizalde@dartmouth.edu}. Partially supported by Simons Collaboration Grant
\#929653
and by Dartmouth-Kalaniyot in the Sciences (DKS) Grant \#1013.
} 
\and Ira M. Gessel~\thanks{Department of Mathematics, Brandeis University, Waltham, MA 02453, USA. E-mail: {\tt gessel@brandeis.edu}.
}
\and 
Yuval Roichman~\thanks{Department of Mathematics, Bar-Ilan University, Ramat-Gan 52900, Israel. 
E-mail: {\tt yuvalr@math.biu.ac.il}.
Partially supported by Dartmouth-Kalaniyot in the Sciences (DKS) Grant \#1013.
}
}

\date{November 11, 2026}

\begin{document}

\maketitle

\begin{abstract}
    We study a variation of the shuffle operation on permutations, which acts on their cycle structure instead of their one-line notation. 
    Given two permutations  $\sigma\in \S_k$ and  $\tau\in \S_{n-k}$, the 
    {\em shuffle-by-cycle} 
    $\shc{\sigma}{\tau}$ is the set of all permutations in $\S_n$ whose set of cycles 
    can be separated into two disjoint subsets: 
    one whose standardization is 
    $\sigma$, and one whose standardization is 
    $\tau$. 
    One can think of the sets $\shc{\sigma}{\tau}$ as refinements of conjugacy classes. 

    We show that if $\sigma$ and $\tau$ have no order-isomorphic cycles, then the descent set has the same distribution, over the set $\shc{\sigma}{\tau}$ of 
    shuffles-by-cycle,
    as over the set $\sigma\shuffle\tau$ of standard shuffles; our proof is bijective.
    This result no longer holds when $\sigma$ and $\tau$ have order-isomorphic cycles. For an arbitrary shuffle-by-cycle, we provide a description    
    of the distribution of the descent set 
    in terms of plethysm of quasisymmetric functions.   
      For   the set $\shc{\sigma}{\sigma}\subset \S_{2n}$, where $\sigma$ is an $n$-cycle,  
    we deduce explicit generating functions for the bivariate distribution of descent number and major index.
\end{abstract}


\tableofcontents

\section{Introduction}

\subsection
{Shuffle-by-cycle}
The main objects of study in this paper are 
shuffles-by-cycle,
which are certain refinements of conjugacy classes in the symmetric group. Let us start with the basic definitions.

Let $\S_n$ denote the symmetric group, that is, the set of all permutations of $[n] := \{1,2,\dots,n\}$. Let $\Co_{(n)}\subset\S_n$ be the set of $n$-cycles, namely, permutations in $\S_n$ that consist of only one cycle.
We will use two different notations for permutations $\pi\in\S_n$: the one-line notation $\pi=\pi_1\pi_2\dots\pi_n$, and the cycle notation, where we write $\pi$ as a product of disjoint cycles of the form $(i,\pi_i,\pi_{\pi_i},\dots)$. 
Given a sequence of distinct positive integers $a_1a_2\dots a_m$, its 
{\em standardization} 
$\st(a_1a_2\dots a_m)\in\S_m$ is obtained by replacing the smallest entry by $1$, the second smallest by $2$, and so on.
We say that two such sequences $a_1a_2\dots a_m$ and $b_1b_2\dots b_m$ are {\em order-isomorphic} if $\st(a_1a_2\dots a_m)=\st(b_1b_2\dots b_m)$, that is, $a_i<a_j$ if and only if $b_i<b_j$ for all $i,j$. 
It is also possible to apply standardization to sets of disjoint cycles, in which case we call it {\em cycle standardization} and denote it by $\st_c$. To apply $\st_c$, we replace the smallest entry in the cycle notation by $1$, the second smallest by $2$, and so on. For example, $\st_c((4,7)(2,5,3))=(3,5)(1,4,2)$.
Two cycles are {\em order-isomorphic} if they have the same cycle standardization.


\begin{definition}\label{def:twist} 
    Let $\sigma\in\S_k$ and $\tau\in\S_{n-k}$. 
    The 
    {\em shuffle-by-cycle}
    of $\sigma$ and $\tau$, denoted by $\shc{\sigma}{\tau}$,
    is the set of all permutations $\pi\in\S_n$ whose set of cycles is the disjoint union of two subsets, 
    one whose cycle standardization equals the cycle form of $\sigma$, and one whose cycle standardization equals the cycle form of $\tau$. 
\end{definition}

\begin{example}\label{ex:A}
    If $\sigma=(1,3)(2)$ and $\tau=(1,4,3)(2,5)$, then $\pi = (\bl{2},\bl{7})(\bl{3})(\rd{1},\rd{6},\rd{5})(\rd{4},\rd{8}) \in \shc{\bl{\sigma}}{\rd{\tau}}$ because $\st_c((\bl{2},\bl{7})(\bl{3}))=\sigma$ and $\st_c((\rd{1},\rd{6},\rd{5})(\rd{4},\rd{8}))=\tau$.
\end{example}

We can also describe $\shc{\sigma}{\tau}$ in terms of the one-line notation, as follows. For $\pi \in \S_n$ and a subset $K\subseteq[n]$, denote by $\pi|_K$ the subsequence of $\pi_1\pi_2\dots\pi_n$ obtained by restricting to the positions in $K$.
Then $\shc{\sigma}{\tau}$ is the set of permutations $\pi\in\S_n$ for which there exists a disjoint union $[n]=K\sqcup L$ such that $\pi(K)=K$, $\pi(L)=L$, $\st(\pi|_K)=\sigma$, and $\st(\pi|_L)=\tau$ in one-line notation. 

\begin{example}
The permutations from Example~\ref{ex:A} in one-line notation are $\sigma=321$ and $\tau=45132$. We have $\pi=\rd{6}\bl{73}\rd{815}\bl{2}\rd{4}$, where $K=\{\bl{2},\bl{3},\bl{7}\}$ and $L=\{\rd{1},\rd{4},\rd{5},\rd{6},\rd{8}\}$.
\end{example}

For comparison, for permutations $\sigma = \sigma_1\sigma_2\dots\sigma_k \in \S_k$ and $\tau = \tau_1\tau_2\dots\tau_{n-k} \in \S_{n-k}$ (written in one-line notation), we denote by $\sigma\shuffle\tau$ the usual shuffle of
$\sigma$ with $\tau'=(\tau_1+k)(\tau_2+k) \dots (\tau_{n-k}+k)$. This is the set of all permutations $\pi\in\S_n$ for which $\sigma$ and $\tau'$ are subsequences of $\pi$. 
Clearly, $|\sigma\shuffle\tau|=\binom{n}{k}$.

\medskip

We first investigate the case when no cycle of $\sigma$ is order-isomorphic to a cycle of $\tau$; we say that $\sigma$ and $\tau$ have {\em no order-isomorphic cycles}. In this case, we have $|\shc{\sigma}{\tau}|=\binom
{n}{k}$, since $\pi\in \shc{\sigma}{\tau}$ is uniquely determined by the partition $[n]=K\sqcup L$.

\begin{example} For $\sigma =(1,3,2) = 312$ and $\tau = (1,2) = 21$ we have
\begin{align*}
\shc{\bl{\sigma}}{\rd{\tau}}&=\{(\bl{1},\bl{3},\bl{2})(\rd{4},\rd{5}),(\bl{1},\bl{4},\bl{2})(\rd{3},\rd{5}),(\bl{1},\bl{5},\bl{2})(\rd{3},\rd{4}),(\bl{1},\bl{4},\bl{3})(\rd{2},\rd{5}),(\bl{1},\bl{5},\bl{3})(\rd{2},\rd{4}),\\ & \qquad (\bl{1},\bl{5},\bl{4})(\rd{2},\rd{3}),(\bl{2},\bl{4},\bl{3})(\rd{1},\rd{5}),(\bl{2},\bl{5},\bl{3})(\rd{1},\rd{4}),(\bl{2},\bl{5},\bl{4})(\rd{1},\rd{3}),(\bl{3},\bl{5},\bl{4})(\rd{1},\rd{2})\}\\
&=\{\bl{312}\rd{54},\bl{41}\rd{5}\bl{2}\rd{3},\bl{51}\rd{43}\bl{2},\bl{4}\rd{5}\bl{13}\rd{2},\bl{5}\rd{4}\bl{1}\rd{2}\bl{3},\bl{5}\rd{32}\bl{14},\rd{5}\bl{423}\rd{1},\rd{4}\bl{52}\rd{1}\bl{3},\rd{3}\bl{5}\rd{1}\bl{24},\rd{21}\bl{534}\},
\end{align*}
$$\bl{\sigma}\shuffle\rd{\tau}=\{\bl{312}\rd{54},\bl{31}\rd{5}\bl{2}\rd{4},\bl{31}\rd{54}\bl{2},\bl{3}\rd{5}\bl{12}\rd{4},\bl{3}\rd{5}\bl{1}\rd{4}\bl{2},\bl{3}\rd{54}\bl{12},\rd{5}\bl{312}\rd{4},\rd{5}\bl{31}\rd{4}\bl{2},\rd{5}\bl{3}\rd{4}\bl{12},\rd{54}\bl{312}\}.$$
\end{example}


We also study the case where  $\sigma\in \S_k$ and $\tau\in \S_{n-k}$ have 
order-isomorphic cycles. 
In this case,  $|\shc{\sigma}{\tau}|<|\sigma\shuffle \tau|$.   
For example, 
if $\sigma$ is the $3$-cycle $(1,3,2)$, then 
\begin{align*}
\shc{\sigma}{\sigma}=\{(1,3,2)(4,6,5),\, (1,4,2)(3,6,5),\, 
(1,5,2)(3,6,4),\,  (1,6,2)(2,5,4),\, (1,4,3)(2,6,5), \\ (1,5,3)(2,6,4),\, 
(1,6,3)(2,5,4),\, (1,5,4)(2,6,3),\, (1,6,4)(2,5,3),\, (1,6,5)(2,4,3)\}. 
\end{align*}



  

    For every $\sigma\in \S_k$ and $\tau\in \S_{n-k}$, 
    denote by $\sigma\oplus\tau \in \S_n$ the permutation whose one-line notation is $\sigma_1 \dots \sigma_k (\tau_1 + k) \dots (\tau_{n-k} + k)$.     
    Denoting the identity permutation in $\S_n$ by $\id_n=12\dots n$, we clearly have
    \[
        \shc{\sigma}{\tau}
        = \{\pi^{-1} (\sigma\oplus\tau) \pi:\ \pi \in \id_k\shuffle \id_{n-k}\}.
    \]
In particular, every 
shuffle-by-cycle
$\shc{\sigma}{\tau}$ is a subset of a conjugacy class. In the other direction, 
every conjugacy class is a disjoint union of suitable 
shuffles-by-cycle;
see Observation~\ref{obs:disjoint_union} below.

\medskip

The 
shuffle-by-cycle
$\shc{\sigma}{\tau}$ was 
initially introduced in~\cite{AR},
in order to extend enumerative results from perfect matchings to all matchings. 
This corresponds to taking $\sigma\in \S_k$ to be a fixed-point-free involution, and $\tau\in \S_{n-k}$ the identity permutation.

\begin{remark} 
    The 
    shuffle-by-cycle
    $\shc{\sigma}{\tau}$ should not be confused with the cyclic shuffle $\sigma\shuffle_{\rm cyc}\tau$ studied in~\cite{AGRR}. 
    For example, $\shc{(1,2,3)}{(1)}=\{(1,2,3)(4),\, (1,2,4)(3),\, (1,3,4)(2),\, (2,3,4)(1)\}$, 
    a set of $3$-cycles 
    in $\S_4$, 
    whereas $(1,2,3)\shuffle_{\rm cyc} (4)=\{(1,2,3,4),\, (1,2,4,3),\, (1,4,2,3)\}$, a set of $4$-cycles in $\S_4$. 
\end{remark}

\subsection{Main results}

The descent set of a permutation $\pi\in\S_n$ is  
\[
    \Des(\pi)
    = \{i: \pi(i) > \pi(i+1)\} \subseteq [n-1].
\]
The fundamental quasisymmetric function corresponding to a subset $J \subseteq [n-1]$ is
\[
    F_{n,J}
    := \sum_{\substack{i_1 \le i_2 \le \cdots \le i_n \\ j\in J \Rightarrow i_j<i_{j+1}}} x_{i_1} x_{i_2} \cdots x_{i_n}.
\]
The set 
$\{F_{n,J} \,:\, J \subseteq [n-1]\}$ 
forms a basis for the quasisymmetric functions of degree $n$;
hence, for any subset $A\subseteq \S_n$,  
the quasisymmetric function $\sum_{\pi \in A} F_{n,\co(\pi)}$  
determines the distribution of the descent set 
on $A$. 

\begin{theorem}\label{thm:main}
    For every $\sigma\in \S_k$ and $\tau\in \S_{n-k}$ with no order-isomorphic cycles, 
    \begin{equation}\label{eq:main}
        \sum_{\pi \in \shc{\sigma}{\tau}} F_{n,\co(\pi)}
        = \sum_{\pi' \in \sigma\shuffle\tau} F_{n, \co(\pi')}.
    \end{equation}
    Equivalently, for any pair $\sigma\in\S_k$ and $\tau\in\S_{n-k}$ with no order-isomorphic cycles, 
    the statistic $\Des$ has the same distribution on the set of
    shuffles-by-cycle
    $\shc{\sigma}{\tau}$ as on the set of standard shuffles $\sigma\shuffle\tau$. 
\end{theorem}

A version of Theorem~\ref{thm:main} was conjectured in~\cite[Conjecture 6.1]{AR}, 
and an outline of an 
algebraic proof was suggested by the third author in~\cite{G}. In this paper we give two detailed proofs, one of which is bijective. Note that the 
natural bijection that maps $(\sigma\oplus\tau)\pi\in \sigma\shuffle \tau$, where $\pi\in \iota_k\shuffle \iota_{n-k}$, to $\pi^{-1}(\sigma\oplus\tau)\pi\in \shc{\sigma}{\tau}$ is not descent-set-preserving. 
Our bijection expands upon and 
refines the Gessel--Reutenauer
bijection \cite{GR} from words to multisets of necklaces. 

\medskip

Recall also the {\em descent number} $\des(\pi):=\card{\Des(\pi)}$ and the {\em major index} $\maj(\pi):=\sum_{i\in \Des(\pi)}i$ of a permutation $\pi\in \S_n$. 
An immediate consequence of Theorem~\ref{thm:main} is that each of the statistics $\des$ and $\maj$ has the same distribution on  $\shc{\sigma}{\tau}$ as on $\sigma\shuffle\tau$. We can then use Stanley's Shuffling Theorem (see, e.g., \cite[Theorem 1.1]{GZ} and references therein) to deduce the following.

\begin{corollary}\label{cor:des-maj} 
    For every $\sigma\in \S_k$ and $\tau\in \S_{n-k}$ with no order-isomorphic cycles, 
    \[
        \sum_{\pi \in \shc{\sigma}{\tau}} q^{\maj(\pi)}
        = q^{\maj(\sigma) + \maj(\tau)}  {{n} \brack {k}}_q
    \]
    and
    \[
        \sum_{\pi \in \shc{\sigma}{\tau}} t^{\des(\pi)}
        = \sum_{i=0}^{n} \binom{k+\des(\tau)-\des(\sigma)}{i-\des(\sigma)} \binom{n-k+\des(\sigma)-\des(\tau)}{i-\des(\tau)}\, t^i .
    \]
\end{corollary}

\medskip

Definition~\ref{def:twist} may be naturally generalized to an arbitrary number of factors.
We need only a restricted version, where all the cycles of each factor are order-isomorphic, and cycles in distinct factors are order-distinct.

\begin{definition}\label{def:multi-shuffle}
    Let $t,n_1,\dots,n_t,k_1,\dots,k_t$ be positive integers, and denote $n := \sum_{i=1}^{t} k_i n_i$.    Let $\sigma_1 \in \S_{n_1},\dots ,\sigma_t \in \S_{n_t}$ be order-distinct (fixed-point-free) cycles.
    The 
    {\em shuffle-by-cycle} $\shuffle_c(\sigma_1^{k_1}, \dots, \sigma_t^{k_t})$     
    is the set of all permutations $\pi\in \S_n$ which have $k_i$ cycles order-isomorphic to $\sigma_i$, for each $1 \le i \le t$.
\end{definition}


For a partition $\nu \vdash n$, let $\Co_\nu$ be the 
conjugacy class in $\S_n$ consisting of all permutations of cycle type $\nu$.  

\begin{observation}\label{obs:disjoint_union}
    Every conjugacy class $\Co_\nu$ 
    in $\S_n$ 
    is a disjoint union of suitable
    $\shuffle_c({\sigma_1^{k_1}, \dots, \sigma_t^{k_t}})$, where $\sigma_1,\dots ,\sigma_t$ are order-distinct cycles.
\end{observation}

\begin{example}
For the partition $\nu = (3,3,2) \vdash 8$,
\[
    \Co_{(3,3,2)}
    = \shuffle_c((1,2,3)^2,(1,2)) \sqcup \shuffle_c((1,2,3),(1,3,2),(1,2)) \sqcup \shuffle_c((1,3,2)^2,(1,2)).
\]
\end{example}

In the next statement, $h_k$ is the complete homogeneous symmetric function of degree~$k$, and $h_k[F_J]$ denotes plethysm of quasisymmetric functions, to be defined in Section~\ref{sec:plethysm} (see also~\cite{MR}). 
The proof of the following general result uses Theorem~\ref{thm:main}. 

\begin{theorem}\label{thm:main2}
For every 
shuffle-by-cycle
$\shuffle_c(\sigma_1^{k_1}, \dots, \sigma_t^{k_t})\subset \S_n$, where  
$\sigma_1,\dots ,\sigma_t$ are order-distinct cycles of lengths $r_1,\dots,r_t$ respectively, 
\[
    \sum_{\pi \in \shuffle_c(\sigma_1^{k_1}, \dots, \sigma_t^{k_t})} F_{n,\co(\pi)}
    = \prod_{i=1}^t h_{k_i}[F_{r_i,\co(\sigma_i)}].
\]    
\end{theorem}

The analogue of Theorem~\ref{thm:main} for $\sigma = \tau$ is false.
To get an analogue of Corollary~\ref{cor:des-maj} for $\sigma = \tau \in \Co_{(n)}$
we apply Theorem~\ref{thm:main2}, 
use a linear operator mapping quasisymmetric functions 
to the ring $\QQ[[t, q]]$ with the Hadamard product in $t$, 
and apply $q$-hypergeometric identities.

\begin{theorem}\label{thm:4}
    For every cycle $\sigma\in \Co_{(n)}$,
    \[
        \sum_{\pi \in \shc{\sigma}{\sigma}} q^{\maj(\pi)} t^{\des(\pi)} 
        = q^{2\maj(\sigma)} t^{\des(\sigma)} \cdot  \sum\limits_{i=0}^n \frac{1}{2}\left({n \brack i}_q^2 + (-1)^i {n \brack i}_{q^2} \right) q^{i^2} t^i. 
    \]
\end{theorem}

We deduce two univariate formulas.

\begin{corollary}\label{t4:3}
    For every cycle $\sigma\in \Co_{(n)}$,
    \[
        \sum_{\pi \in \shc{\sigma}{\sigma}} q^{\maj(\pi)}
        = q^{2\maj(\sigma)} \cdot \frac{1}{2} \left( {2n \brack n}_q+  \prod\limits_{i=1}^n (1-q^{2i-1})
        \right).
    \]
\end{corollary}

\begin{corollary}\label{cor4:2}
    For every cycle $\sigma\in \Co_{(n)}$,
    \[
        \sum_{\pi \in \shc{\sigma}{\sigma}} t^{\des(\pi)} 
        = t^{\des(\sigma)} \cdot \sum_{i=0}^n a_{i,n} t^i,
    \]
    where
    \[
        a_{i,n}
        = \begin{cases}
            \displaystyle\binom{\binom{n}{i}+1}{2}, & \ \text{if $i$ is even};\\
            {}\\
            \displaystyle\binom{\binom{n}{i}}{2}, & \ \text{if $i$ is odd}. 
        \end{cases}
    \]
\end{corollary}



\section{The case of
no order-isomorphic cycles}

\subsection{
The Gessel--Reutenauer bijection from words to necklaces}

In this subsection we recall the   bijection of Gessel and Reutenauer~\cite{GR} from words on a totally ordered alphabet to multisets
of necklaces. 



Let $w=w_1,\dots,w_n$ be a word in a totally ordered alphabet $X=\{x_1<x_2<\dots\}$. 
The standard permutation of $w$ is the unique permutation
$\sigma=\st(w)\in \S_n$ such that 
\[
\sigma_i<\sigma_j \Longleftrightarrow w_i< w_j \ \text{or}\ \left( w_i=w_j \ \text{and}\ i<j\right).  
\]
See Example~\ref{example:21} below. 

Two words $u, v$ in the totally ordered alphabet $X$ are conjugate if for some words $s, t$ one has $u = st$ and $v = ts$. A word w is primitive if it 
is not a proper power of another word. A necklace is the conjugacy class of a primitive word. 
A multiset of necklaces is a mapping $M$ from the set of necklaces into $\NN$, with
finite support. 
The cycle structure of a multiset $M$ of necklaces is the partition where, for each $i$, the number of parts equal to $i$ is the number of necklaces in $M$ of length $i$.   
The evaluation of a necklace $v$ is the monomial $\ev(v):=x_1^{n_1}x_2^{n_2}\cdots \in \ZZ[X]$,
where $n_i$ is the number of occurrences of the variable $x_i$ in $v$; the evaluation
of a multiset $M$ of necklaces is the product of the evaluations of the
necklaces in $M$.


In \cite[Lemma~3.4]{GR}, Gessel and Reutenauer describe an evaluation-preserving bijection $U$ between words and multisets of necklaces, such that, for any word $w$, the cycle structure of $\st(w)$ and $U(w)$ are equal. In some of the examples below it will be convenient to use the alphabet $\{\a<\b<\c<\dots\}$ instead of $\{x_1<x_2<\dots\}$.

Let $w$ be a word, let $M=U(w)$ be the corresponding multiset of necklaces, and let $\pi=\st(w)^{-1}$. We will use the notation $\perm(M)=\pi$, and note that $M$ and $\pi$ have the same cycle structure.

Here is a description of the bijection $U$. Given a word $w=w_1w_2\dots w_n$, let $\pi=\st(w)^{-1}$ (in~\cite{GR}, the bijection is described using $\st(w)$, but for us it will be more convenient to use its inverse).
Consider the cycle form of $\pi$, and replace each element $j$ with $w_{\pi(j)}$ (note that, when writing $\pi^{-1}$ in one-line notation, the position of element $j$ is the same as the position of $w_{\pi(j)}$ in $w$).
Let $U(w)$ be the resulting multiset of necklaces.

The inverse map $U^{-1}$ can be described as follows. Given a multiset of necklaces $M$, label each element of each necklace by the periodic sequence obtained by reading the necklace starting at that element. Then, order these sequences lexicographically (if some necklaces appear with multiplicity, break the ties in a consistent way by first ordering the repeated necklaces). This order assigns a number from $1$ to $n$ to each element of each necklace, which produces the permutation $\pi=\perm(M)$ in cycle form. To recover $w$, write $\pi^{-1}$ in one-line notation and replace each entry by the corresponding letter in $M$.

\begin{example}\label{example:21}
    Let $w=\d\c\a\c\a\c\b\a$. We first find its standard permutation:
    \begin{equation}\label{eq:st}\begin{array}{rccccccccl}
        w&=&\d&\c&\a&\c&\a&\c&\b&\a  \\
        \st(w)&=&8&5&1&6&2&7&4&3,
    \end{array}\end{equation}
    whose inverse $\pi=\st(w)^{-1}=35872461$ has cycle form  $(4,7,6)(2,5)(1,3,8)$. Replacing each element with the corresponding letter from Equation~\eqref{eq:st}, we get the multiset of necklaces $U(w)=(\b,\c,\c)(\a,\c)(\a,\a,\d)$.

    To compute the inverse map, we label each element of $M=(\b,\c,\c)(\a,\c)(\a,\a,\d)$ by the periodic sequence starting at that element:
    $$(\b\c\dots,\c\c\dots,\c\b\dots)(\a\c\dots,\c\a\dots)(\a\a\dots,\a\d\dots,\d\a\dots).$$
    Ordering these sequences lexicographically gives
    $\pi=\perm(M)=(4,7,6)(2,5)(1,3,8)$. Writing $\pi^{-1}=85162743$ and replacing each entry by the corresponding letter in $M$, we recover $w=\d\c\a\c\a\c\b\a$.
\end{example}

\subsection{Proofs of Theorem~\ref{thm:main} and Corollary~\ref{cor:des-maj}}\label{sec:proof_1.4}

Next we state a property of the map $\perm$ that will be useful in the proofs. We use $\sqcup$ to denote the disjoint union of multisets, where  the multiplicities of each element are added.

\begin{lemma}\label{lem:2multisets}
    Let $\sigma\in\S_k$ and $\tau\in\S_{n-k}$, and let $M$ be a multiset of necklaces. Then $\perm(M)\in \shc{\sigma}{\tau}$ if and only if there exist multisets of necklaces $M_1$ and $M_2$ such that $M=M_1\sqcup M_2$, $\perm(M_1)=\sigma$ and $\perm(M_2)=\tau$.

    Additionally, if $\sigma$ and $\tau$ have no order-isomorphic cycles, then the decomposition $M=M_1\sqcup M_2$ is unique.
\end{lemma}

\begin{proof}
    To obtain the permutation $\perm(M)$ associated to a multiset of necklaces $M$, we label each element of each necklace by the periodic sequence starting at that element, and order these sequences lexicographically. This gives the cycle form of $\perm(M)$. The labels of the elements of a particular necklace, as well as the relative order among these labels, only depend on that necklace, and not on the rest of $M$.

    Suppose that $M=M_1\sqcup M_2$ with $\perm(M_1)=\sigma$, and $\perm(M_2)=\tau$. Then the relative order of the labels assigned to $M_1$ is the same as in the computation of $\perm(M_1)=\sigma$, and similarly for $M_2$. Thus, the cycle form of $\pi$ can be separated into the set of cycles (coming from $M_1$) whose entries are in the same relative order as the cycle form of $\sigma$, and the set of cycles (coming from $M_2$) whose entries are in the same relative order as the cycle form of $\tau$. In other words, $\pi\in \shc{\sigma}{\tau}$.

    Conversely, suppose now that $M$ is a multiset of necklaces such that $\perm(M)\in \shc{\sigma}{\tau}$. By definition of $\shc{\sigma}{\tau}$, the cycle form of $\perm(M)$ can be separated into two sets of cycles, one whose standardization is $\sigma$, and one whose standardization is $\tau$. Let $M_1$ consist of those necklaces in $M$ corresponding to the cycles in the first set, and let $M_2$ consist of the remaining necklaces, so that $M=M_1\sqcup M_2$. By the first paragraph, $\perm(M_1)=\sigma$ and $\perm(M_2)=\tau$.

    If $\sigma$ and $\tau$ have no order-isomorphic cycles, then for each necklace $N$ in $M$, the cycle $\perm(N)$ is order-isomorphic to a cycle of $\sigma$ or $\tau$, but not both. This uniquely determines the decomposition $M=M_1\sqcup M_2$.
\end{proof}

Another ingredient in the proof will be the fact, shown in~\cite[Exer.\ 7.93]{EC2}, that for any $\sigma\in\S_k$ and $\tau\in\S_{n-k}$,
\begin{equation}\label{eq:F-shuffle}
\sum_{\pi\in\sigma\shuffle\tau}F_{n,\co(\pi)}=F_{k,\co(\sigma)}F_{n-k,\co(\tau)}.
\end{equation}

\begin{proof}[First proof of Theorem~\ref{thm:main}]
    By \cite[Lem.~3.2]{GR} and \cite[Lem.~3.4]{GR},   for any $\pi\in \S_n$ we have
    \[    F_{n,\co(\pi)}=\sum_{w:\st(w)=\pi^{-1}}\ev(w)=\sum_{M:\perm(M)=\pi}\ev(M).
    \]
    Thus, we can write the right-hand side of Equation~\eqref{eq:main} as $$\sum_{\pi'\in \shc{\sigma}{\tau}}F_{n,\co(\pi')}=\sum_{M:\perm(M)\in \shc{\sigma}{\tau}}\ev(M).$$

    By Lemma~\ref{lem:2multisets}, multisets of necklaces $M$ such that $\perm(M)\in \shc{\sigma}{\tau}$ are those that can be decomposed as $M=M_1\sqcup M_2$ with $\perm(M_1)=\sigma$ and $\perm(M_1)=\tau$. Since $\sigma$ and $\tau$ have no order-isomorphic cycles, this decomposition is unique. 
    Noting that $\ev(M)=\ev(M_1)\ev(M_2)$, it follows that     
    \begin{align*}
        \sum_{\pi' \in \shc{\sigma}{\tau}} F_{n,\co(\pi')}
        &= \sum_{(M_1,M_2) \,:\, \perm(M_1) = \sigma,\, \perm(M_2) = \tau} \ev(M_1) \ev(M_2)\\
        &= \sum_{M_1 \,:\, \perm(M_1) = \sigma} \ev(M_1) \sum_{M_2 \,:\, \perm(M_2) = \tau} \ev(M_2) \\
        &= F_{k,\co(\sigma)} F_{n-k,\co(\tau)}.
    \end{align*}

Finally, by Equation~\eqref{eq:F-shuffle}, this product also equals the right-hand of Equation~\eqref{eq:main}.
\end{proof}

The above proof 
can be turned into a bijection, combining both Gessel and Reutenauer's bijection $U$ and its inverse $U^{-1}$.

\begin{proof}[Second proof of Theorem~\ref{thm:main}]
    By~\cite[Lemma~3.2]{GR}, we can write Equation~\eqref{eq:main} as
    $$\sum_{w:\st(w)^{-1}\in\sigma\shuffle\tau}\ev(w)
    =\sum_{w':\st(w')^{-1}\in \shc{\sigma}{\tau}}\ev(w').$$
    We will describe an evaluation-preserving bijection
    \begin{equation}\label{eq:ww'}
    \bij:\{w:\st(w)^{-1}\in\sigma\shuffle\tau\}\to\{w':\st(w')^{-1}\in \shc{\sigma}{\tau}\}.\end{equation}
    Afterwards we will also show that the map $\bij$ can be easily translated to a bijection for permutations.

    Given a word $w=w_1w_2\dots w_n$ in the left-hand side of Equation~\eqref{eq:ww'}, let $u=w_1w_2\dots w_k$ and $v=w_{k+1}w_{k+2}\dots w_n$. By the definition of shuffles and standardization, we have $\st(u)^{-1}=\sigma$ and $\st(v)^{-1}=\tau$. Applying Gessel and Reutenauer's bijection $U$ to each of $u$ and $v$, and noting that $\sigma$ and $\tau$ are cycles, we get two multisets of necklaces $U(u)$ and $U(v)$. Now let $M=U(u)\sqcup U(v)$, and apply the inverse of $U$ to get a word $w'=U^{-1}(M)$. By Lemma~\ref{lem:2multisets}, we have $\st(w')^{-1}=\perm(M)\in \shc{\sigma}{\tau}$, so $w'$ belongs to the right-hand side of Equation~\eqref{eq:ww'}. Let $\bij(w)=w'$.

    The inverse of the map $\bij^{-1}$ can be described by reversing each of the above steps. Given $w'$ in the right-hand side of Equation~\eqref{eq:ww'}, let $\pi'=\st(w')^{-1}\in \shc{\sigma}{\tau}$. The multiset of necklaces $M=U(w')$ satisfies $\perm(M)=\pi'$, so by Lemma~\ref{lem:2multisets}, it can be uniquely decomposed as $M=M_1\sqcup M_2$ with $\perm(M_1)=\sigma$ and $\perm(M_2)=\tau$. 
    Applying $U^{-1}$ to each of these multisets, we obtain words $u=U^{-1}(M_1)$ and $v=U^{-1}(M_2)$ such that $\st(u)^{-1}=\sigma$ and $\st(v)^{-1}=\tau$. It follows that the concatenation $w=uv$ satisfies $\st(w)^{-1}\in\sigma\shuffle\tau$. We have $\bij^{-1}(w')=w$.
    \medskip

    A variation of the above map $\bij$ gives a bijection in terms of permutations. Specifically, for any fixed subset $S\subseteq[n-1]$, we will describe a bijection
     $$\bijS:\{\pi\in\sigma\shuffle\tau:\Des(\pi)\subseteq S\}
    \to\{\pi'\in \shc{\sigma}{\tau}:\Des(\pi')\subseteq S\},$$
    from where it follows that the statistic $\Des$ has the same distribution on $\sigma\shuffle\tau$ as on $\shc{\sigma}{\tau}$.
    We write $S=\{s_1,s_2,\dots,s_{m-1}\}$, where $1\le s_1<s_2<\dots<s_{m-1}\le n-1$. 

    Given $\pi\in\sigma\shuffle\tau$ with $\Des(\pi)\subseteq S$, first construct a word $w$ by writing $\pi^{-1}$ in one-line notation and replacing the $s_1$ smallest entries with $x_1$, the next $s_2-s_1$ smallest entries with $x_2$, the next $s_3-s_2$ smallest entries with $x_3$, and so on, until finally replacing the $n-s_{m-1}$ largest entries with $x_m$. By construction, $\pi=\st(w)^{-1}$. As in the description of $\bij$, let $u=w_1w_2\dots w_k$ and $v=w_{k+1}w_{k+2}\dots w_n$, and let $M=U(u)\sqcup U(v)$. Now define $\bijS(\pi)=\perm(M)$.
    
    To describe the inverse map, start with $\pi'\in \shc{\sigma}{\tau}$, written in cycle form, and replace the $s_1$ smallest entries with $x_1$, the next $s_2-s_1$ smallest entries with $x_2$, the next $s_3-s_2$ smallest entries with $x_3$, and so on, to get multiset of necklaces $M$ which can be uniquely decomposed as $M=M_1\sqcup M_2$ with $\perm(M_1)=\sigma$ and $\perm(M_2)=\tau$.
    Let $u=U^{-1}(M_1)$, $v=U^{-1}(M_2)$, and $w=uv$. Then $\bijS^{-1}(\pi)=\st(w)^{-1}$.    
\end{proof}

Below are some examples of $\bijS$ and $\bijS^{-1}$.

\begin{example}
    Let $\sigma=(1,3,2,4)=3421\in\S_4$ and $\tau=(1,3,2)=312\in\S_3$, and let $S=\{2,3,6\}$. We will apply $\bijS$ to $\pi=3742561\in\sigma\shuffle\tau$, which has $\Des(\pi)=S$.
    We take $\pi^{-1}=7413562$ and replace entries $1,2$ with $\a$, entry $3$ with $\b$, entries $4,5,6$ with $\c$, and entry $7$ with $\d$. This gives $w=\d\c\a\b\c\c\a$, and so $u=\d\c\a\b$ and $v=\c\c\a$, which satisfy $\st(u)^{-1}=\sigma$ and $\st(v)^{-1}=\tau$. Making the substitutions $j\mapsto u_{\sigma_j}$ in the cycle notation of $\sigma$, we get the necklace $U(u)=(\a,\c,\b,\d)$, and substituting $j\mapsto v_{\tau_j}$ in the cycle notation of $\tau$ gives $U(v)=(\a,\c,\c)$. Labeling each element in the multiset of necklaces $M=(\a,\c,\b,\d)(\a,\c,\c)$ with the periodic sequence starting at that element, we get
    $$(\a\c\b\dots,\c\b\d\dots,\b\d\a\dots,\d\a\c\dots)(\a\c\c\dots,\c\c\a\dots,\c\a\c\dots),$$
    and ordering the elements lexicographically yields $\bijS(\pi)=\perm(M)=(1,5,3,7)(2,6,4)=5672341$, which has descent set $\{3,6\}\subseteq S$.
\end{example}

\begin{example}
    Let $\sigma$ and $\tau$ be as in the previous example, and let now $S=\{1,2,4,5,6\}$. Let $\pi'=(2,5,4,7)(1,6,3)=6517432\in \shc{\sigma}{\tau}$, which has $\Des(\pi')=S$. To apply $\bijS^{-1}$ to it, we first make the replacements $1\mapsto\a$, $2\mapsto\b$, $\{3,4\}\mapsto\c$, $5\mapsto\d$, $6\mapsto\e$, $7\mapsto\f$ in the cycle form of $\pi'$, which yields the pair of necklaces
    $M_1=(\b,\d,\c,\f)$ and $M_2=(\a,\e,\c)$. To apply $U^{-1}$ to $M_1$, since we already know that $\perm(M_1)=\sigma=(1,3,2,4)$, we simply write $\sigma^{-1}=4312$ and replace each element with the corresponding letter from $M_1$, which gives $u=\f\d\b\c$. 
    Similarly, to apply $U^{-1}$ to $M_2$, which satisfies
    $\perm(M_2)=\tau=(1,3,2)$, we write $\tau^{-1}=231$ and replace each element with the corresponding letter from $M_2$, which gives $v=\c\e\a$. Now $w=uv=\f\d\b\c\c\e\a$, whose standardization is $\st(w)=7523461$, so $\bijS^{-1}(\pi')=\st(w)^{-1}=7345261$, which has descent set $\{1,4,6\}\subseteq S$.
\end{example}

\begin{example}
    Let $\sigma=(2,3,5)(1,4,6)=435621\in\S_6$ and $\tau=(3)(1,4,2)=4132\in\S_4$, and let $S=\{1,3,6,8,9\}$. We will apply $\bijS$ to $\pi=\bl{4\,3}\,\rd{10}\,\bl{5\,6}\,\rd{7}\,\bl{2}\,\rd{9\,8}\,\bl{1}\in\bl{\sigma}\shuffle\rd{\tau}$, which has $\Des(\pi)=S$.
    We take $\pi^{-1}=10\,7\,2\,1\,4\,5\,6\,9\,8\,3$ and make the replacements $1\mapsto\a$, $\{2,3\}\mapsto\b$, $\{4,5,6\}\mapsto\c$, $\{7,8\}\mapsto\d$, $9\mapsto\e$, $10\mapsto\f$. This gives $w=\f\d\b\a\c\c\c\e\d\b$, and so $u=\f\d\b\a\c\c$ and $v=\c\e\d\b$, which satisfy $\st(u)^{-1}=\sigma$ and $\st(v)^{-1}=\tau$. Making the substitutions $j\mapsto u_{\sigma_j}$ in the cycle notation of $\sigma$, and the substitutions $j\mapsto v_{\tau_j}$ in the cycle notation of $\tau$, we get the multisets of necklaces $U(u)=(\b,\c,\d)(\a,\c,\f)$ and $U(v)=(\d)(\b,\e,\c)$. Labeling each element in $M=U(u)\sqcup U(v)=(\b,\c,\d)(\a,\c,\f)(\d)(\b,\e,\c)$ with the periodic sequence starting at that element, we get
    $$(\b\c\dots,\c\d\dots,\d\b\dots)(\a\c\dots,\c\f\dots,\f\a\dots)(\d\d\dots)(\b\e\dots,\e\c\dots,\c\b\dots),$$
    and ordering the elements lexicographically yields 
    $$\bijS(\pi)=\perm(M)=(\bl{2},\bl{5},\bl{7})(\bl{1},\bl{6},\bl{10})(\rd{8})(\rd{3},\rd{9},\rd{4})=\bl{6\,5}\,\rd{9\,3}\,\bl{7\,10\,2}\,\rd{8\,4}\,\bl{1},$$ which has descent set $\{1,3,6,8,9\}=S$.
\end{example}

\begin{proof}[Proof of Corollary~\ref{cor:des-maj}]
The formula for $\maj$
follows from 
Theorem~\ref{thm:main} together with~\cite[Theorem 3.1]{GG}. 
The formula for $\des$
follows from 
Theorem~\ref{thm:main} together with~\cite[Corollary 5.7]{AGRR}. 
\end{proof}

\section{The general 
case}

\subsection{Plethysm of quasisymmetric functions}\label{sec:plethysm}

Plethysm for quasisymmetric functions was defined  in~\cite{MR}; see also~\cite[Def.\ A2.6]{EC2} for the symmetric case.
In particular, if $f$ is a symmetric function and $g = \sum_i u_i$ is any formal power series which is a sum of monomials $u_i$ with no numerical coefficients, say of the form $x_1^{j_1}x_2^{j_2}\cdots$ or $q^j$, then 
\[
    f[g]
    := f(u_1, u_2,\dots),
\]
that is, $f[g]$ is obtained by replacing the variables of $f$ by the monomials $u_i$.

\subsection{Proof of Theorem~\ref{thm:main2}}


    Recall, from Definition~\ref{def:multi-shuffle}, that for a cycle $\sigma \in \S_r$, the shuffle-by-cycle
    $\shuffle_c(\sigma^k)\subset \S_{rk}$ is the set of all permutations $\pi\in \S_{rk}$ which have $k$ cycles order-isomorphic to $\sigma$.
    
    \begin{lemma}\label{lem:hF} 
    Let $\sigma \in \Co_{(r)}$ be an $r$-cycle in $\S_r$ with descent set $\Des(\sigma)$. 
    The quasisymmetric descent set generating function for the
    shuffle-by-cycle
    $\shuffle_c(\sigma^k)\subset \S_{rk}$ 
    is equal to the plethysm 
    \[
        \sum_{\pi \in \shuffle_c(\sigma^k)} F_{rk,\Des(\pi)}
        = h_k[F_{r,\Des(\sigma)}],
    \]
    where $h_{k}$ is the $k$-th complete homogeneous symmetric function.
\end{lemma}

\begin{proof} 
    As in the first proof of Theorem~\ref{thm:main} in Subsection~\ref{sec:proof_1.4} above, 
    \[
        F_{r,\co(\sigma)}
        = \sum_{M:\,\perm(M)=\sigma} \ev(M). 
    \]
    Similarly,
    \[
        \sum_{\pi \in \shuffle_c(\sigma^k)} F_{rk,\co(\pi)}
        = \sum_{\pi \in \shuffle_c(\sigma^k)} \sum_{M:\,\perm(M)=\pi} \ev(M).
    \]
    This double summation is actually over all the multisets $M$ consisting of $k$ (not necessarily distinct) necklaces, each mapped to $\sigma$ under $\perm$. 
    It follows that 
    \[
        \sum_{\pi \in \shuffle_c(\sigma^k)} F_{rk,\co(\pi)}
        = h_k[F_{r,\co(\sigma)}].
        \qedhere
    \]
\end{proof}

Now we are ready to prove  Theorem~\ref{thm:main2}.

\begin{proof}[Proof of Theorem~\ref{thm:main2}]
    By Theorem~\ref{thm:main}, if $\sigma\in\S_k$ and $\tau\in\S_{n-k}$ have no order-isomorphic cycles, then 
    \[
        \sum_{\pi \in \shc{\sigma}{\tau}} F_{n,\co(\pi)}
        = \sum_{\pi' \in \sigma \shuffle \tau} F_{n,\co(\pi')}
         = F_{k,\co(\sigma)} F_{n-k,\co(\tau)},
    \]
    where the last equality follows from Equation~\eqref{eq:F-shuffle}.
    
    Thus, by iteration, if $\sigma_1, \dots ,\sigma_t$ are order-distinct cycles of lengths $r_1,\dots, r_t$ then 
    \[
        \sum_{\pi \in \shuffle_c(\sigma_1^{k_1}, \dots, \sigma_t^{k_t})} F_{n,\co(\pi)}
        = \prod_{i=1}^t \Biggl(\sum_{\pi \in \shuffle_c(\sigma_i^{k_i})} F_{r_i k_i, \co(\pi)} \Biggr). 
    \]
    Lemma~\ref{lem:hF} completes the proof.
\end{proof}

\section{The 
case of equal cycles}

In this section we prove Theorem~\ref{thm:4}. 

First, note that 
Theorem~\ref{thm:main} does not hold when $\sigma$ and $\tau$ have order-isomorphic cycles. For the case of $\sigma = \tau \in \Co_k$, we instead have the following 
expression.

\begin{lemma}\label{thm:sigma=tau}
    For every cycle $\sigma\in\Co_{(n)}$, 
\[
        \sum_{\pi \in \shc{\sigma}{\sigma}} F_{2n,\co(\pi)}=
         \frac{1}{2} \left(F_{n,\Des(\sigma)}^2 + \left. F_{n,\co(\sigma)} \right|_{x_i\mapsto x_i^2} \right)
        = \frac{1}{2} \left( \sum_{\pi' \in \sigma\shuffle\sigma} F_{2n,\co(\pi')} + \left. F_{n,\co(\sigma)} \right|_{x_i\mapsto x_i^2} \right) ,
\]
    where $\left. f \right|_{x_i\mapsto x_i^2} := f(x_1^2, x_2^2, \ldots)$.

    
    
 
\end{lemma}





\begin{proof}[First proof of Lemma~\ref{thm:sigma=tau}]
    Recall that 
    \[
        h_2
        = \frac{1}{2} \left( p_1^2 + p_2 \right),
    \]
    where $p_k$ is the $k$-th power sum symmetric function. 
    By Lemma~\ref{lem:hF}, 
    \[
        \sum_{\pi \in \shc{\sigma}{\sigma}} F_{2n,\co(\pi)}
        = h_2 \left[F_{n,\Des(\sigma)} \right]
        = \frac{1}{2} \left( p_1^2 \left[ F_{n,\Des(\sigma)} \right] + p_2 \left[ F_{n,\Des(\sigma)} \right] \right)
        = \frac{1}{2} \left( F_{n,\Des(\sigma)}^2 + \left. F_{n,\co(\sigma)} \right|_{x_i\mapsto x_i^2} \right). 
    \]
    Finally, by Equation~\eqref{eq:F-shuffle},
    \[  
        F_{n,\Des(\sigma)}^2
        = \sum_{\pi \in \sigma\shuffle\sigma} F_{2n,\Des(\pi)}.
        \qedhere
    \]
\end{proof}

\begin{proof}[Second proof of Lemma~\ref{thm:sigma=tau}]
    We proceed as in the first proof of Theorem~\ref{thm:main}. In this case, each of $M_1$ and $M_2$ consists of one necklace.
    However, the decomposition $M=M_1\sqcup M_2$ is no longer unique, and we have to be careful not to double count pairs of distinct necklaces, since they correspond to the same multiset $M$. Separating the cases when $M_1=M_2$ and $M_1\neq M_2$, we get
    \begin{align*}
    \sum_{\pi\in \shc{\sigma}{\sigma}}F_{2n,\co(\pi)}
    &=\sum_{N:\perm(N)=\sigma}\ev(N)^2
    +\frac{1}{2}\sum_{\substack{(M_1,M_2):M_1\neq M_2\\ \perm(M_1)=\perm(M_2)=\sigma}}\ev(M_1)\ev(M_2)\\
    &=\frac{1}{2}\sum_{N:\perm(N)=\sigma}\ev(N)^2
    +\frac{1}{2}\sum_{\substack{(M_1,M_2):\\ \perm(M_1)=\perm(M_2)=\sigma}}\ev(M_1)\ev(M_2)\\
    &=\frac{1}{2}\left(\left.F_{n,\co(\sigma)}\right|_{x_i\mapsto x_i^2}+F_{n,\co(\sigma)}^2\right).
    \qedhere
    \end{align*}
\end{proof}



 \bigskip





Consider the linear operator $\Psi$ on quasisymmetric functions defined by
\[
    \Psi(f(x_1,x_2,\ldots))
    := \sum_{i=0}^\infty f(1,q,\dots,q^{i-1},0,\ldots) t^i.
\]
Recall the notation $(t;q)_{r} := (1-t)(1-tq)\cdots(1-tq^{r-1})$ and $(q)_{r} := (q;q)_{r}$, for $r \ge 1$. 

The following lemma reformulates~\cite[Theorem 4.5(c)]{GZ}. 
We give a self-contained proof.

\begin{lemma}\label{lem:psiF}
For any subset $J=\{j_1,\ldots,j_k\} \subseteq [n-1]$,
\[
    \Psi(F_{n,J})
    = \frac{q^{\comaj(J)} t^{|J|+1}}{(t;q)_{n+1}}
\]
where $\comaj(J) := \sum_{i=1}^{k} (n-j_i) = \sum_{s=1}^{n} |J \cap [s-1]|$.     
\end{lemma}


\begin{proof} 
    We follow the arguments in the proof of~\cite[Lemma 7.19.10]{EC2}. By definition,
    \[
        F_{n,J}(1,q,\ldots,q^{i-1},0,\ldots)
        = \sum_{\substack{1 \le i_1 \le \cdots \le i_n \le i \\ i_j < i_{j+1} \text{ if } j \in J}} q^{(i_1-1) + \cdots + (i_n-1)} .
    \]
    Denoting
    \[
        r_s := (i_s-1) - |J \cap [s-1]|
    \qquad (1 \le s \le n),
    \]
    the conditions $1 \le i_1 \le \cdots \le i_n \le i$ and $i_j < i_{j+1}$ for $j \in J$ are equivalent to $0 \le r_1 \le \dots \le r_n \le (i-1) -|J|$. 
    Thus, by~\cite[Proposition 1.7.3]{EC1},
    \begin{align*}
        F_{n,J}(1,q,\ldots,q^{i-1},0,\ldots)
        &= \sum_{0 \le r_1 \le \cdots \le r_n \le i-1-|J|} \prod_{s=1}^{n} q^{r_s + |J \cap [s-1]|} \\
        &= q^{\comaj(J)}
        \cdot \sum_{0 \le r_1 \le \cdots \le r_n \le i-1-|J|} q^{r_1 + \cdots + r_n}. \\ 
        &= q^{\comaj(J)} \cdot {n + i-1-|J| \brack n}_q .
    \end{align*}
    Therefore, by the formula between Equations (1.86) and (1.87) in~\cite{EC1},
    \[
        \Psi(F_{n,J})
        = q^{\comaj(J)} \cdot \sum_{i=|J|+1}^{\infty} {n + i-1-|J| \brack n}_q t^i 
        = \frac{q^{\comaj(J)} t^{|J|+1}}{(t;q)_{n+1}}.
        \qedhere
    \]
\end{proof}



%

\begin{proof}[Proof of Theorem~\ref{thm:4}] 
    For an $n$-cycle $\sigma\in \Co_{(n)}$ denote
    \[
        R_2(\sigma)
        := \sum_{\pi \in \shc{\sigma}{\sigma}} q^{\comaj(\pi)} t^{\des(\pi)+1}.
    \]

    By Lemmas~\ref{lem:psiF} and~\ref{thm:sigma=tau}, 
    \begin{align*}        
         R_2(\sigma)
         &= \sum_{\pi \in \shc{\sigma}{\sigma}} (t;q)_{2n+1} \Psi(F_{2n,\co(\pi)})
        = \frac{(t;q)_{2n+1}}{2}         \left( \Psi(F_{n,\Des(\sigma)}^2) + \Psi \left( \!\left. F_{n,\co(\sigma)} \right|_{x_i\mapsto x_i^2} \right) \right). 
    \end{align*}
    Now
    \[
        \Psi \left( \!\left. F_{n,\co(\sigma)} \right|_{x_i\mapsto x_i^2} \right)
        = \left. \Psi(F_{n,\co(\sigma)}) \right|_{q \mapsto q^2}
        = \frac{q^{2\comaj(\sigma)} t^{\des(\sigma)+1}}{(t;q^2)_{n+1}}.
    \]
    Also, by the proof of Lemma~\ref{lem:psiF},
    \[
        F_{n,J}(1,q,\ldots,q^{i-1},0,\ldots)
        = q^{\comaj(J)} \cdot {n + i-1-|J| \brack n}_q .
    \]
    Therefore
    \begin{align*}
        \Psi(F_{n,\Des(\sigma)}^2)
        &= \sum_{i=0}^\infty F_{n,\Des(\sigma)}^2(1,q,\dots,q^{i-1},0,\ldots) t^i
        = \sum_{i = \des(\sigma)+1}^{\infty} q^{2\comaj(\sigma)} {n + i-1-\des(\sigma) \brack n}_q^2 t^i \\
        &= q^{2\comaj(\sigma)} t^{\des(\sigma)+1} \cdot \sum_{j = 0}^{\infty}  {n + j \brack n}_q^2 t^j         
    \end{align*}
    and 
    \begin{equation}\label{eq:R2}
        R_2(\sigma)
        = \frac{1}{2} q^{2\comaj(\sigma)} t^{\des(\sigma)+1} (t;q)_{2n+1} \cdot \left( \frac{1}{(t;q^2)_{n+1}} + \sum_{j = 0}^{\infty}  {n + j \brack n}_q^2 t^j \right) .    
    \end{equation}
    We shall use $q$-hypergeometric identities to evaluate~\eqref{eq:R2} as a polynomial in $t$.
    
    First of all, by the $q$-binomial theorem~\cite[Equation (1.87)]{EC1}
    \begin{equation}\label{eq:R2a}
        \frac{(t;q)_{2n+1}}{(t;q^2)_{n+1}}
        = \prod_{i=1}^{n} (1 - tq^{2i-1})
        = (tq;q^2)_{n}
        = \sum_{i=0}^{n} {n \brack i}_{q^2} (q^2)^{\binom{i}{2}} (-tq)^i
        = \sum_{i=0}^{n} {n \brack i}_{q^2} q^{i^2} (-t)^i .
    \end{equation}
    On the other hand, substituting $a = b = q^{n+1}$ and $c = q$ in Heine's $q$-analogue of Euler's transformation formula~\cite[Equation (1.4.3)]{GaR} 
    \[
        \sum_{j=0}^{\infty} \frac{(a;q)_j (b;q)_j}{(c;q)_j (q;q)_j} z^j
        = \frac{(abz/c;q)_\infty}{(z;q)_\infty} \sum_{i=0}^{\infty} \frac{(c/a;q)_i (c/b;q)_i}{(c;q)_i (q;q)_i} (abz/c)^i
    \]    
    yields
    \[
        \sum_{j=0}^{\infty} \frac{(q^{n+1};q)_j^2}{(q;q)_j^2} z^j
        = \frac{(q^{2n+1}z;q)_\infty}{(z;q)_\infty} \sum_{i=0}^{\infty} \frac{(q^{-n};q)_i^2}{(q;q)_i^2} (q^{2n+1}z)^i .
    \]
    Since $(q^{-n};q)_i = 0$ for $i > n$, this reduces to 
    \[
        \sum_{j=0}^{\infty} {n + j \brack n}_q^2 z^j
        = \frac{1}{(z;q)_{2n+1}} \sum_{i=0}^{n} \frac{(q^{-n};q)_i^2}{(q;q)_i^2} (q^{2n+1}z)^i .
    \]
    Using
    \[
        (q^{-n};q)_i
        = (-1)^i q^{-i(2n+1-i)/2} (q^{n+1-i};q)_i
    \]
    gives
    \begin{equation}\label{eq:R2b}
        (z;q)_{2n+1} \cdot \sum_{j=0}^{\infty} {n + j \brack n}_q^2 z^j
        = \sum_{i=0}^{n} \frac{q^{-i(2n+1-i)} (q^{n+1-i};q)_i^2}{(q;q)_i^2} (q^{2n+1}z)^i 
        = \sum_{i=0}^{n} q^{i^2} {n \brack i}_q^2 z^i .
    \end{equation}
    Substituting~\eqref{eq:R2a} and~\eqref{eq:R2b} (with $z = t$) into~\eqref{eq:R2} gives
    \[
        R_2(\sigma)
        = \frac{1}{2} q^{2\comaj(\sigma)} t^{\des(\sigma)+1} \cdot \sum_{i=0}^{n} \left( {n \brack i}_{q^2} q^{i^2} (-t)^i + {n \brack i}_q^2 q^{i^2} t^i \right) .    
    \]
    To get the major index instead of the comajor index, replace $\sigma$ by $w_0 \sigma w_0$, where $w_0 = [n, n-1, \dots, 1] \in \S_n$, and similarly for $\pi$. 
\end{proof}

\begin{proof}[Proof of Corollary~\ref{t4:3}.]
Letting $t=1$ in Theorem~\ref{thm:4} yields, 
for every $n$-cycle $\sigma\in \Co_{(n)}$,
\[
\sum\limits_{\pi\in \shc{\sigma}{\sigma}} q^{\maj(\pi)}= q^{2\maj(\sigma)}\cdot  \sum\limits_{i=0}^n \frac{q^{i^2}}{2}\left({n \brack i}_q^2 + (-1)^i {n \brack i}_{q^2} \right).  
\]
By the $q$-Chu-Vandermode identity~\cite[Exercise~1.100]{EC1}, 
\[
\sum\limits_{i=0}^n q^{i^2} {n \brack i}_q^2={2n \brack n}_q. 
\]
By the $q$-binomial theorem,
\[
\sum_{i=0}^n (-1)^i q^{i^2} {n \brack i}_{q^2}=\prod\limits_{i=1}^r (1-q^{2i-1}).
\]
This completes the proof.
\end{proof}

\begin{proof}[Proof of Corollary~\ref{cor4:2}]
Setting $q=1$ in Theorem~\ref{thm:4} yields, 
for every $n$-cycle $\sigma\in \Co_{(n)}$,    
\[
    \sum_{\pi \in \shc{\sigma}{\sigma}} t^{\des(\pi)} 
    = \frac{t^{\des(\sigma)}}{2} \cdot \sum_{i=0}^n \left( {n \choose i}^2 + (-1)^i {n \choose i} \right) t^i
    = t^{\des(\sigma)} \cdot \sum_{0}^n a_{i,n} t^i,
\]
where $a_{i,n}$ is equal to $\binom{\binom{n}{i}}{2}$ if $i$ is odd, and to $\binom{\binom{n}{i}+1}{2}$ if $i$ is even.
\end{proof}

\section{Final remarks 
}\label{sec:final}



It is possible to generalize Theorem~\ref{thm:4} and get closed formulas for $k>2$, but they will be more complicated. 

Recall from~\cite{GZ} the operation of {\em Hadamard product $*$} on formal power series in $t$, defined by  
\[
    \sum_{i=0}^\infty a_i t^i * \sum_{i=0}^\infty b_i t^i 
    := \sum_{i=0}^\infty a_i b_i t^i. 
\]
We consider the algebra of formal power series in $t$ and $q$, in which the multiplication $*$ is the usual one in $q$, 
but 
in $t$ it is the Hadamard product. 
For example,
\[
    (1 + q t + q^2 t^2) * (1 + q^2 t + q^4 t^2) 
    = 1 + q^3 t + q^6 t^2.
\]
Generalization of the proof of Theorem~\ref{thm:4} to any $k\ge 2$  involves Hadamard products of  power series.  
For other applications of this algebra to permutation enumeration see~\cite{GZ}.


\begin{proposition}
    For every $n$-cycle $\sigma \in \Co_{(n)}$ and $k>0$
    \begin{equation}\label{eq:511}
        \sum_{\pi \in \shuffle_c(\sigma^k)} q^{\maj(\pi)} t^{\des(\pi)} 
        =  q^{k\maj(\sigma)} t^{\des(\sigma)} (t;q)_{nk+1} h_k \left[  \frac{1}{(t;q)_{n+1}} 
        \right], 
    \end{equation}
    where $h_k \left[ \frac{1}{(t;q)_{n+1}} \right]$ is computed in the ring $\QQ[[t, q]]$, with the Hadamard product in $t$. 
    Explicitly,
    \[
        h_k \left[ \frac{1}{(t;q)_{n+1}} \right] 
        = \sum_{\lambda \vdash k} \frac{1}{z_\lambda} \sum_{j=0}^\infty t^j \prod_{i=1}^{\infty} {{n+j} \brack n}_{q^{i}}^{m_i}, 
    \]
    where the outer summation is over all partitions $\lambda 
    = (1^{m_1}2^{m_2}\cdots)$ of $k$, and $z_\lambda := \prod_i m_i! i^{m_i}$.
\end{proposition}

\begin{proof}
Let   $\sigma\in \Co_{(n)}$ be an $n$-cycle in $\S_n$.
By Lemma~\ref{lem:psiF} and Lemma~\ref{lem:hF}, 
\[
    \sum_{\pi \in \shuffle_c(\sigma^k)} q^{\comaj(\pi)} t^{\des(\pi)+1} 
    = (t;q)_{nk+1} \cdot \sum_{\pi \in \shuffle_c({\sigma^k})} \Psi(F_{nk,\co(\pi)})
    = (t;q)_{nk+1} \cdot \Psi(h_k[F_{n,\co(\sigma)}]).
\]
Clearly, $\Psi$ commutes with plethysm in the ring $\QQ[[t, q]]$, with the Hadamard product in $t$. 
Thus, by Lemma~\ref{lem:psiF}, 
\[
    \Psi(h_k[F_{n,\co(\sigma)}]) 
    = h_k[\Psi(F_{n,\co(\sigma)})] 
    = h_k \left[ \frac{q^{\comaj(\sigma)} t^{\des(\sigma)+1}}{(t;q)_{n+1}} \right]
    = q^{k\comaj(\sigma)} t^{\des(\sigma)+1} h_k \left[ \frac{1}{(t;q)_{n+1}} \right]. 
\]
To get the major index instead of the comajor index, replace $\sigma$ by $w_0 \sigma w_0$. This completes the proof of Equation~\eqref{eq:511}.

To get an explicit expression for $h_k \left[ \frac{1}{(t;q)_{n+1}} \right]$, note that by the formula between Equations (1.86) and (1.87) in~\cite{EC1}, 
\[
\frac{1}{(t;q)_{n+1}}=\sum\limits_{j=0}^\infty {n+j \brack n}_q t^j. 
\]
Recall also that $h_k = \sum_{\lambda \vdash k} \frac{p_\lambda}{z_\lambda}$ and that, for $\lambda = (1^{m_1} 2^{m_2} \cdots)$,
\[
    p_\lambda \left[ {{n+j}\brack n}_q t^j \right]
    = \prod_{i=1}^\infty {{n+j}\brack n}_{q^i}^{m_i} t^j.  
\]
This completes the proof.
\end{proof}

A natural problem is to find analogues of Theorem~\ref{thm:4} for other permutation statistics. 
    We believe that methods similar to those used in this paper may be useful to 
    other permutation statistics determined  by descents, such as the number of peaks. 

\medskip

Finally, for an application of Theorem~\ref{thm:main} to the enumeration of matchings, see~\cite[Theorem 6.6]{AR}. 
More applications of this kind are desired.

\end{document}